\documentclass[12pt]{amsart}

\usepackage{epsfig}
\usepackage{amscd}
\usepackage[mathscr]{eucal}
\usepackage{amssymb}
\usepackage{amsxtra}
\usepackage{amsmath}
\usepackage[all]{xy}

\theoremstyle{plain}

\theoremstyle{definition}

\theoremstyle{remark}

\newcommand{\cntrs}{\setcounter{thm}{0}
  \renewcommand{\thethm}{\thesection.\Alph{thm}}}

\begin{document}

\title{A Note on Lagrangian Vanishing Sphere Bundles}
\author{Yochay Jerby} 

%\bibliographystyle{plain}

%----------------------------------------------------------------------
%
% Abstract
%
%\begin{abstract}
% ...
%\end{abstract}

\maketitle

%----------------------------------------------------------------------
%
% Beginning of text
%

\section{Introduction} \cntrs
\label{S:Intro}

 \hspace{-0.6cm} The study of Lagrangian submanifolds (and the restrictions imposed on them) plays a pivotal
  role in Symplectic
 topology. For this reason it is good to have an abundance of
 examples of Lagrangian submanifolds. Generating such examples,
 however, is not always a simple task as, naturally, each such example
 involves a unique construction.

 \hspace{-0.6cm} A classical way to construct a Lagrangian in a symplectic manifold $\Sigma$, which in a sense
 is considered as part of the "folklore" of the field, is
  to let $\Sigma$ appear as a smooth fiber in
   a Lefschetz fibration. If this is possible the singularities of
  the fibration induce Lagrangian spheres in $\Sigma$ and these
  spheres, in turn, are representatives of the corresponding vanishing cycles
  in the homology of $\Sigma$
 , see for instance
  \cite{Ar:monodromy,Se}.

  \hspace{-0.6cm} In this short work
our aim is twofold: The first would be describing a generalization
of the above mentioned construction to the "Morse-Bott" case. This
would lead, whenever such a degeneration exists, to the existence of Lagrangian sphere bundles rather
than just spheres. The second part of the paper would be to study
the topological restrictions on such Lagrangian sphere bundles
implied by the theory of Floer homology for Lagrangian
intersections and to illustrate the techniques involved. 

\hspace{-0.6cm} Let $X$ be a K$\ddot{\textrm{a}}$hler manifold and
let $\Omega$ be a K$\ddot{\textrm{a}}$hler form on $X$ which we
will assume from now on to be fixed. Our setting in this paper is
the following:

\bigskip

 \hspace{-0.6cm} \bf Definition: \rm A proper holomorphic map
 $\pi: X \rightarrow \mathbb{D}$ map from a K$\ddot{\textrm{a}}$hler 
 manifold $X$ to the unit disc $ \mathbb{D} \subset \mathbb{C}$
 is called a Morse-Bott degeneration of a K$\ddot{\textrm{a}}$hler manifold 
 $\Sigma$ 
 if:

 \bigskip

 (1) The only critical value of $\pi$ is $ 0 \in \mathbb{D}$ and Crit($\pi$) is a
 submanifold of $X$.

 \bigskip

 (2) The holomorphic Hessian of $\pi$ is a
   non-degenerate quadratic form when

   \hspace{0.5cm} restricted to the normal bundle of Crit($\pi$) in X.
  
  \bigskip
  
  (3) $\Sigma$ is biholomorphic to a nonsingular fiber of $\pi$. 
   
\bigskip

\hspace{-0.6cm}  Note that if $Crit(\pi)$ is a point the above
definition reduces to give a description of the local model of a
classical Lefschetz fibration around a singularity. Our first
result is hence the following:

\bigskip

\hspace{-0.6cm} \bf Theorem A: \rm  Let $\pi : X
\rightarrow \mathbb{D}$ be a Morse-Bott degeneration and let $L
\subset Crit( \pi)$ be a compact Lagrangian submanifold. Every
smooth fiber $\Sigma$ of $\pi$ contains a
Lagrangian submanifold $N(L)$ diffeomorphic to a sphere bundle
over $L$.

\bigskip

\hspace{-0.6cm}  We shall refer to the Lagrangian submanifold
$N(L)$ as the Lagrangian "vanishing necklace" over $L$ induced by
the degeneration $\pi$. First, we have the following
 straightforward example of a Morse-Bott degeneration: 

\bigskip

\hspace{-0.6cm} \bf Example: \rm  Let $\pi : X \rightarrow \mathbb{D}$ be a Morse-Bott degeneration of $\Sigma$ with an isolated singularity (i.e a Lefschetz fibration)
 and let $p : E \rightarrow X$ be a holomorphic projective bundle on $X$ of rank $k$. Consider 
the map $\widetilde{ \pi}_E : E \rightarrow \mathbb{D}$ given by $\widetilde{\pi}_E(x)=\pi(p(x))$. Then $\widetilde{\pi}_E$ is a Morse-Bott degeneration with $Crit(\widetilde{\pi}_E) = \mathbb{P}^k$. 

\bigskip

\hspace{-0.6cm}  From now on let us consider the case of a
Morse-Bott degeneration with $Crit(\pi)=\mathbb{P}^k$. Note that the fiber of the Morse-Bott degeneration 
with $Crit(\pi)=\mathbb{P}^k$ described above is by itself a projective bundle. It seems, however, that if we require $h^{1,1}(\Sigma)=1$
 an example of such a degeneration is increasingly hard to find. Indeed, in the second part of the paper we will apply
 methods of symplectic topology, particularly the theory of Floer homology of Lagrangian submanifolds, to show that Morse-Bott 
 degenerations with $Crit(\pi)=\mathbb{P}^k$ admit non-trivial restrictions. 

\hspace{-0.6cm} We shall first apply Theorem A to show that under the assumption $Crit(\pi)=\mathbb{P}^k$  one can find 
a K$\ddot{\textrm{a}}$hler manifold $(A, \omega)$ such that $dim_{\mathbb{C}}A = k+1$ and $ \Sigma \times A$ admits a displaceable 
Lagrangian necklace $N(L)$ which is topologically a sphere bundle over $S^{2k+1}$. The
construction of this Lagrangian necklace is described in detail in
section 4. One main feature of the Lagrangian $N(L)$ is that due
to the affine factor it would be hamiltoninaly displaceable from
itself. In particular, by the fundamental properties of Floer homology, whenever $N(L)$ admits 
a well defined Floer homology algebra $HF(N(L),N(L))$, it would inevitably be trivial. 

\hspace{-0.6cm} Recall that a K$\ddot{\textrm{a}}$hler manifold is said to be Fano if $c_1(\Sigma)$, the first Chern class of $\Sigma$, 
can be represented by a K$\ddot{\textrm{a}}$hler form. One way to ensure that the Floer homology algebra $HF(N(L),N(L))$ is well defined is to
require the of monotonicity of $N(L)$, see \cite{Oh:HF1,Oh:spectral} which would hold if we require $\Sigma$ to be a Fano manifold with 
$h^{1,1}(\Sigma)=1$. Assume that this is indeed the case and define the number $$C_{\Sigma}:= min \left \{ \left  < c_1(\Sigma) , A \right > \mid A \in \pi_2(\Sigma) \textrm{ such that 
} 
0< \left< c_1(\Sigma) , A \right > \right \} $$ to which we refer as the minimal Chern number of $\Sigma$. We obtain the following restrictions:

\bigskip

\hspace{-0.6cm} \bf Theorem B: \rm Let $\Sigma$ be a Fano manifold of $dim_{\mathbb{C}}\Sigma=n$ with $h^{1,1}(\Sigma)=1$ and minimal Chern 
number $C_{\Sigma}$. Suppose that $\Sigma$ appears as a fibre of a Morse-Bott degeneration $\pi : X \rightarrow \mathbb{D}$ with 
$Crit(\pi) \simeq \mathbb{P}^k$ and $k \neq 0 , n-1$. The following restrictions hold:  

\bigskip

\hspace{-0.6cm} (1) If $ n \neq 3k+1$ then one of the following holds: 

\bigskip

(a) $C_{\Sigma} \mid k+1$. 

\bigskip

(b) $2 C_{\Sigma} \mid n-k+1$. 

\bigskip

(c) $ 2 C_{\Sigma} \mid n+k+2$. Furthermore $C_{\Sigma} \mid 2k+1 $ if $n > 3k+1$ and $C_{\Sigma} \mid 2k+2$ if $n < 3k+1$. 

\bigskip

\hspace{-0.6cm} (2) If $n=3k+1$ then $ C_{\Sigma} \mid n-k+1$. 

\bigskip

\hspace{-0.6cm} It is interesting to note that we have the following example of a Morse-Bott degeneration with 
$Crit(\pi)=\mathbb{P}^{n-1}$, which is the case not covered by our theorem:

\bigskip

\hspace{-0.6cm} \bf Example: \rm  Let $\Sigma$ be a K$\ddot{\textrm{a}}$hler manifold of dimension $dim_{\mathbb{C}} \Sigma = n$ and let $Y \subset \Sigma$ 
be a submanifold. Consider the manifold $X =\Sigma \times \mathbb{D}$ and let $p_Y: \widetilde{X}_Y \rightarrow X$ be the blow up of $X$ along the $Y \times \left \{ 0 \right \}$. Denote by $E$ be the exceptional 
divisor of $\widetilde{X}_Y$. The map 
$\pi : \widetilde{X}_Y \rightarrow \mathbb{D}$ given by $\pi = pr_{\mathbb{D}} \circ p_Y$ is a Morse-Bott degeneration of $\Sigma$ with $Crit(\pi)=E$. In particular. if $Y$ 
is a point $Crit(\pi) \simeq \mathbb{P}^{n-1}$.    

\bigskip

\hspace{-0.6cm} In view of Theorem B, it is tempting to suggest that it is impossible to include a Fano manifold $\Sigma$ with $h^{1,1}(\Sigma)=1$
as a fibre of a Morse degeneration with $Crit(\pi)=\mathbb{P}^k$. Finally, let us conclude with the following remark: 

\bigskip

\hspace{-0.6cm} \bf Remark \rm (Degenerations and pencils): In projective geometry,
degenerations with isolated singularities exist in abundance as they arise naturally
from Lefschetz pencils which could be associated to any projective
manifold $X$. It is interesting to note that such a
straightforward correspondence between degenerations and pencils does
not seem to hold in the non-isolated, Morse-Bott case, due to what
appears to be essential reasons. This would be discussed in detail
at the end of the paper.

\bigskip

\hspace{-0.6cm} The rest of the paper is organized as follows: In
section 2 we prove Theorem A. Section 3 is devoted to an overview
of required results on Floer homology of monotone Lagrangian
submanifolds. In section 4 we discuss Morse-Bott degenerations
with $Crit(\pi)=\mathbb{P}^k$ and prove
Theorem B. We discuss the obstructions to obtain Morse-Bott degenerations from 
pencils in section 5.

\bigskip

\hspace{-0.6cm} \emph{Acknowledgements:} The author would like to
thank his advisor, Professor Paul Biran, for his support and
guidance during the preparation of this work.

\section{Proof of the Lagrangian Necklace Theorem} \label{Sb:Imp}

\hspace{-0.6cm} This section is devoted to the proof of the
following theorem:

\bigskip

\hspace{-0.6cm} \bf Theorem 2.1: \rm  Let $\pi : X \rightarrow
\mathbb{D}$ be a Morse-Bott degeneration and let $L \subset Crit(
\pi)$ be a Lagrangian submanifold. Every smooth fiber
$\Sigma_z=\pi^{-1}(z)$ for $0 \neq z \in \mathbb{D}$ contains a
Lagrangian submanifold $N_z(L)$ diffeomorphic to a sphere bundle
over $L$.

\bigskip

\hspace{-0.6cm} We follow the lines of Donaldson's proof for the Morse case which
appears in \cite{Se}. We refer the reader to \cite{Au-Br} for relevant facts on
Morse-Bott theory. Set $F=Re(\pi)$ and $H=Im(\pi)$ and denote by
$J$ the complex structure of $X$. Since $\pi$ is holomorphic the
negative gradient flow of $F$ with respect to the metric $g_J(
\cdot , \cdot) = \Omega( \cdot , J \cdot)$ is the same as the
Hamiltonian flow of $H$ with respect to the
K$\ddot{\textrm{a}}$hler form $\Omega$. Denote this flow by $
 \phi_t $ and let $S(L)$ be the stable submanifold of $L$ under this flow, see \cite{Kat:Has}.

\bigskip

\hspace{-0.6cm} \bf Lemma 2.2: \rm $S(L)$ is a lagrangian
submanifold in $(X, \Omega)$.

\bigskip

\hspace{-0.6cm} \bf Proof: \rm Let $ \Theta: S(L) \rightarrow L$
be the end point map given by $ \Theta(x): = lim_{t \rightarrow
 \infty} \phi_t(x)$. In a small neighborhood $U \subset S(L)$ of $L$ the
 map $ \Theta$ is a locally trivial fibration over $L$,
see \cite{Kat:Has}. Let us show that $U$ is isotropic. Indeed, for every $ x
\in U$ and $v_1,v_2 \in T_x(U)$ we have
$$ \Omega_x(v_1,v_2) = \Omega_{\pi(x)}(d \Theta(v_1), d
\Theta(v_2))=0$$ because $ \phi_t$ is
Hamiltonian with respect to $ \Omega$ and $L$ is Lagrangian. Moreover, the
fact that $U$ is isotropic implies that $S(L)$ is isotropic since
$\phi_t$ is Hamiltonian and every point $x \in S(L)$ is sent to
$U$ by $ \phi_t$ for $t$ large enouge.

\hspace{-0.6cm} Finally recall that
$$dim_{\mathbb{R}}S(L)=dim_{\mathbb{R}}L+index(F)$$ Where $index(F)$ is the Morse-Bott index, see [3].
A simple computation shows that
$$ indexF= \frac{1}{2} \left ( dim_{\mathbb{R}} X-dim_{\mathbb{R}} Crit(\pi) \right ) =
\frac{1}{2}dim_{\mathbb{R}} X-dim_{\mathbb{R}}L$$ Thus
$dim_{\mathbb{R}}S(L)=\frac{1}{2}dim_{\mathbb{R}} X$ and hence $S(L)$ is
Lagrangian. $\square$

\bigskip

\hspace{-0.6cm} \emph{Continuation of the proof of Theorem 2.1:} Denote $N_{\epsilon}(L) = F^{-1}(\epsilon) \cap
S(L)$ for $ \epsilon>0$ small enough. By Morse-Bott theory
$N_{\epsilon}(L)$ is a sphere bundle over $L$, of dimension
$dim_{\mathbb{R}}N_{\epsilon}(L)=dim_{\mathbb{R}}S(L)-1$. As $\phi_t$ is the Hamiltonian
flow of $H$, the function $H$ is constant along the flow lines of
$\phi_t$. We deduce that $S(L) \subset H^{-1}(0)$. Consequently
$$N_{\epsilon}(L) \subset F^{-1}(\epsilon) \cap H^{-1}(0)=
\pi^{-1}(\epsilon)=\Sigma_{\epsilon}$$

\hspace{-0.6cm} This proves that the fibre $\Sigma_{\epsilon}$ contains
the Lagrangian submanifold $N_{\epsilon}(L)$ which is a sphere
bundle over $L$. Finally, by Moser argument all fibers $\Sigma_z=
\pi^{-1}(z)$, $z \neq 0$ , are symplectomorphic, hence they all
contain such Lagrangians. $\square$

\bigskip

\hspace{-0.6cm} \bf Definition 2.2: \rm Let $\pi : X \rightarrow
\mathbb{D}$ be a Morse-Bott degeneration and $L \subset Crit(
\pi)$ be a Lagrangian submanifold. We refer to $N_z(L) \subset
\Sigma_z$ as the Lagrangian necklace over $L$ in the fibre $\Sigma_z$.

\bigskip

\hspace{-0.6cm} \bf Remark: \rm In case $Crit(\pi)$ is a finite
number of points and $L$ coincides with this set of points then
$N_z(L)$ is just the set of Lagrangian spheres representing the
vanishing cycles as appears in \cite{Ar:monodromy,Se}.

\section{Brief Overview of Floer and Quantum Homology for Monotone Lagrangian Submanifolds} \label{Sb:Imp}

\hspace{-0.6cm} An essential ingredient in our approach is the theory of Floer homology for Lagrangian submanifolds.
This theory was originally introduced by Floer, for some class of Lagrangian submanifolds, in order to study problems of Lagrangian intersections, see e.g \cite{Fl:Morse-theory}.
The theory was later extended by Oh \cite{Oh:HF1} for the class of so called monotone Lagrangian submanifolds. We shall not go here into the details of the this homology theory,
but only indicate our algebraic conventions as well as some essential properties of Floer homology. A detailed account on the subject can be found in
\cite{Fl:Morse-theory,Oh:HF1}, see also \cite{Bi:Lag,Bi:ICM2005,Bi:LQH}.

\hspace{-0.6cm} Let $(M, \omega)$ be a closed symplectic manifold and $L \subset M$ a monotone closed Lagrangian submanifold. Monotone
means that for all $\alpha \in \pi_2(M,L)$ we have $$ \omega(\alpha) >0 \Leftrightarrow  \mu(\alpha)>0$$ Where $\mu : \pi_2 (M,L) \rightarrow \mathbb{Z}$ is the Maslov class homomorphism. We denote by $$N_L=min \left \{ \mu(\alpha) \mid \mu(\alpha)>0, \alpha \in \pi_2(M,L) \right \}$$ The minimal Maslov index of $L$, so that $image(\mu) \subset N_L \cdot \mathbb{Z}$. We shall assume
from now on that $L$ is monotone and $N_L \geq 2$. Under these assumptions one can associate to $L$ an important invariant called the
Floer homology $HF_{\ast}(L,L)$. This is a graded vector space over the field $\mathbb{Z}_2$. It is the homology of a certain chain complex
defined by counting solutions to a (non-linear) elliptic PDE with boundary conditions on L (called pseudo-holomorphic curves). The resulting homology
measures the obstruction to Hamiltonianly displace $L$, i.e for the existence of $\varphi \in Ham(M, \omega)$ s.t $\varphi(L) \cap L = \phi$.

\hspace{-0.6cm} Below we shall use essentially Oh's version of Floer homology as developed in \cite{Oh:HF1} only that our $HF_{\ast}(L,L)$ will be $\mathbb{Z}$-graded
as in  \cite{Bi:Lag,Bi:ICM2005,Bi:LQH} and $N_L$-periodic, i.e $HF_{\ast+N_L}(L,L)  =  HF_{\ast}(L,L)$.

\hspace{-0.6cm} Denote by $\Lambda=\mathbb{Z}_2 [t^{-1},t]$ the ring of Laurent polynomials in $t$ with coefficents in $\mathbb{Z}_2$. We grade $\Lambda$
by declaring the degree of $t$ to be $\mid t \mid =- N_L$. Then our $HF_{\ast}(L,L)$ becomes a graded module over $\Lambda$, and we have
$$ t \cdot HF_i(L,L) = HF_{i-N_L} (L,L) $$ For all $i \in \mathbb{Z}$.

\hspace{-0.6cm} Call a Lagrangian submanifold $L \subset M$ \em{displaceable} \rm if there exists a Hamitonian diffemorphism $\varphi \in Ham(M, \omega)$
such that $\varphi(L) \cap L = \phi$. In what follows we shall use the following important vanishing property of $HF$ : \sl{If $L$ is displaceable then $HF(L,L)=0$.}

\hspace{-0.6cm} \rm Another important ingredient in our Floer-theoretical considerations is a spectral sequence which is useful for computing $HF$. This
sequence was introduced by Oh \cite{Oh:spectral}, but here we shall follow a variant of it as described by Biran in \cite{Bi:Lag}. The main properties of this sequence that
are relevant for our applications are listed in the following theorem.

\bigskip

\hspace{-0.6cm} \bf Theorem \rm (see Oh \cite{Oh:spectral}, Biran \cite{Bi:Lag}): There exists a spectral sequence $\left
\{ E_{p,q}^r,d^r \right \}_{r \geq 1}$ with the following properties:

\bigskip

\hspace{-0.6cm} (1) $E_{p,q}^1 \simeq H_{p+q-pN_L}(L ; \mathbb{Z}_2)   \cdot t^{-p}$ for every $p,q$.

\bigskip

\hspace{-0.6cm} (2) $E^r_{p+1,q} \simeq E^r_{p,q+1-N_L} t^{-1} $ for every $p,q,r$.

\bigskip

\hspace{-0.6cm} (3) $\left \{ E_{p,q}^r, d^r \right \} $ collapses after a finite number of pages
and converges to $HF_{\ast}(L,L)$.

\bigskip

\hspace{-0.6cm} In particualr, if $L$ is displaceable then $E_{p,q}^{\infty}=0$ for every $p,q$.

\section{Morse-Bott Degenerations with Critical Locus $\mathbb{P}^k$} \label{Sb:Imp}

\hspace{-0.6cm} Let $\Sigma$ be a K$\ddot{\textrm{a}}$hler manifold of $dim_{\mathbb{C}} \Sigma =n$. In this section we consider the special case of Morse-Bott degenerations $\pi : X \rightarrow \mathbb{D}$ of $\Sigma$ with $Crit(\pi) \simeq \mathbb{P}^k$. Let $ B^{2n}(r) \subset \mathbb{C}^n$ be the standard symplectic ball 
of radius $r>0$ with symplectic form $\omega_0$ induced by the standard symplectic form on $\mathbb{C}^n$ 
and let $S^{2n-1}(r) = \partial B^{2n}(r)$ be the corresponding sphere. In \cite{ALP} it was observed that, for any $r>1$, the 
  manifold $\mathbb{P}^k \times
B^{2k+2}(r)$ contains the embedded Lagrangian sphere
$$ L = \left \{ (h(z), \bar{z}) \mid z \in S^{2k+1}(1) \right \}
\subset (\mathbb{P}^k \times B^{2k+2}(r), \omega_{FS} \oplus
\omega_0)
$$ where $h : \mathbb{C}^{k+1} \setminus \left \{ 0 \right \}
\rightarrow \mathbb{P}^k$ is the Hopf map given by $z \mapsto
\mathbb{C} \cdot z$ and $\omega_{FS}$ is the Fubini-Study form on $\mathbb{P}^k$ normalized such that $\int_{\mathbb{P}^1} \omega_{FS} = \pi$. 

\hspace{-0.6cm} Let $(A, \omega)$ be a compact closed symplectic manifold of $dim_{\mathbb{C}}A=k+1$ such that there exists 
a symplectic embedding of the standard symplectic ball $(B^{2k+2}(r), \omega_0)$ into $A$ for $r>>0$. Consider the stabilization $$\widetilde{\pi} :
X \times
A \rightarrow \mathbb{D} $$ of the map
of $\pi$ given by $(z,u) \mapsto \pi(z)$. This
induces a Morse-Bott degeneration of $\Sigma \times A$ with $Crit(\widetilde{\pi}) =
\mathbb{P}^k \times A$. Thus, for any K$\ddot{\textrm{a}}$hler form
$\Omega$ on $X$ such that $\Omega
\mid_{\mathbb{P}^k} = \omega_{FS}$ we get a Lagrangian sphere $L \subset Crit(\widetilde{\pi})$ and by Theorem 2.1 we have a 
Lagrangian Necklace $N(L) \subset (\Sigma \times
A , \Omega \mid_{ \Sigma} \times \omega )$. We have:

\bigskip

\hspace{-0.6cm} \bf Proposition 4.1: \rm The Lagrangian necklace $N(L) \subset \Sigma \times A $ is a Hamiltonianly displaceable, simply connected, Lagrangian 
sub-manifold which is topologically a sphere bundle over $S^{2k+1}$. 

\bigskip

\hspace{-0.6cm} \bf Proof: \rm Note that the Lagrangian sphere $L \subset \mathbb{P}^k \times B^{2k+2}(4)$ can be displaced from itself by Hamiltonian isotopy $\varphi$ which 
is the identity on the $\mathbb{P}^k$ factor and acting as a Hamitonian isotopy in the $B^{2k+2}(4)$ factor which displaces $B^{2k+2}(1) \subset B^{2k+2}(4)$ from itself.  
Denote by $L'=\varphi(L)$. In particular, for $\epsilon>0$ small enough, $N_{\epsilon}(L) \cap N_{\epsilon}(L') = \phi$ since $L \cap L' = \phi$.  $\square$  

\bigskip

\hspace{-0.6cm} We further have: 

\bigskip

\hspace{-0.6cm} \bf Lemma 4.3: \rm Let $\Sigma$ be a Fano manifold with $dim_{\mathbb{C}} \Sigma=n$ which satisfies $h^{1,1}(\Sigma)=1$. Assume $\Sigma$ appears as a fibre in a degeneration $\pi : X \rightarrow 
\mathbb{D}$ with $Crit(\pi)=\mathbb{P}^k$. $1 \leq k \leq n-2 $ then $$HF(N(L),N(L))=0$$

\bigskip

\hspace{-0.6cm} \bf Proof: \rm The condition $h^{1,1}(\Sigma)=1$ assures that the form on $\Sigma$ is monotone. In particular, since under the above conditions $N(L)$ 
is simply connected, it would be a monotone Lagrangian submanifold of $\Sigma$. Hence $N(L)$ has a well defined Floer homology. By Proposition 4.1 $N(L)$ is 
Hamiltonianly displaceable and thus $HF(N(L),N(L))=0$. $\square$ 

\bigskip

\hspace{-0.6cm} We are now in position to prove the following theorem: 

\bigskip

\hspace{-0.6cm} \bf Theorem 4.4: \rm Let $\Sigma$ be a Fano manifold of $dim_{\mathbb{C}}\Sigma=n$ with $h^{1,1}(\Sigma)=1$ and minimal Chern 
number $C_{\Sigma}$. Suppose that $\Sigma$ appears as a fibre of a Morse-Bott degeneration $\pi : X \rightarrow \mathbb{D}$ with 
$Crit(\pi) \simeq \mathbb{P}^k$ and $k \neq 0 , n-1$. The following restrictions hold:  

\bigskip

\hspace{-0.6cm} (1) If $ n \neq 3k+1$ then one of the following holds: 

\bigskip

(a) $C_{\Sigma} \mid k+1$. 

\bigskip

(b) $2 C_{\Sigma} \mid n-k+1$. 

\bigskip

(c) $ 2 C_{\Sigma} \mid n+k+2$. Furthermore $C_{\Sigma} \mid 2k+1 $ if $n > 3k+1$ and $C_{\Sigma} \mid 2k+2$ if $n < 3k+1$. 

\bigskip

\hspace{-0.6cm} (2) If $n=3k+1$ then $ C_{\Sigma} \mid n-k+1$. 

\bigskip

\hspace{-0.6cm} \bf Proof: \rm First, by standard considerations, the Maslov index of $N(L)$ satisfies $N_L = 2 C_{\Sigma}$. By the Gysin sequence for sphere bundles, 
the singular homology $H_i(N(L),\mathbb{Z}_2)=0$ can only be nonzero if $$i = 0,2k+1,n-k,n+k+1 $$ Let us consider the spectral sequence described in section 3. 
Taking $p=q=0$, the differential of the spectral sequence has the following form $$...\rightarrow E^r_{r,-r+1} \rightarrow E^r_{0,0} \rightarrow E^r_{-r,r-1} \rightarrow ... $$ 
For $r \geq 0$ Moreover, for degree reasons we have $E^r_{r,-r+1}=0$ for all $r \geq 0$. 
 Thus, since $E_{0,0}^1 \simeq \mathbb{Z}_2$, in order for $E^{r+1}_{0,0}$ to vanish for some $r \geq 1$ we must have $E^r_{-r,r-1} \neq 0$. But this 
can only happen in one of the following cases: 

\bigskip

\hspace{-0.6cm} (i) $2rC_{\Sigma}-1=2k+1$ for some $r \geq 0$. 

\bigskip

\hspace{-0.6cm} (ii) $2rC_{\Sigma}-1=n-k$ for some $r \geq 0$. 

\bigskip

\hspace{-0.6cm} (iii) 

\bigskip

(iii.1) $2r_1C_{\Sigma}-1=n+k+1$ for some $r \geq 0$.

\bigskip

(iii.2.1) $(2k+1)+2r_2C_{\Sigma}-1 =n-k$ if $3k+1<n$.

\bigskip

(iii.2.2) $(n-k)+2r_2C_{\Sigma}-1 =2k+1$ if $n<3k+1$.

\bigskip

\hspace{-0.6cm} Finally, note that (iii) can hold only if $n \neq 3k+1$ and in the complementary case (i) and (ii) coincide. This proves (1) and (2). $\square$

\bigskip

\section{On The Construction of Morse-Bott Degenerations} \label{Sb:Imp}

\hspace{-0.6cm} A common way to obtain degenerations with isolated singular points is
 by considering Lefschetz fibrations which in turn arise from Lefschetz pencils. In the non-isolated case, however, 
  there seems to be no such straightforward transition from pencils to fibrations. In this section we would like to explain this difference between 
  the isolated and non-isolated cases. We refer the reader to \cite{Vo} for a complete treatment on the geometry of Lefscehtz pencils.
    
  \hspace{-0.6cm} Let $X \subset \mathbb{P}^N$ be an algebraic manifold and denote by $(\mathbb{P}^N)^{\ast}$ the dual projective space parametrizing hyperplanes 
  in $ \mathbb{P}^N$. To $X$ one can associate the discriminant variety given by $$ X^{\ast} = \left \{ H \mid \Sigma_H \textrm{ is singular} \right \}  \subset (\mathbb{P}^N)^{\ast} $$ 
  where  $\Sigma_H=X \cap H $ is the hyperplane section corresponding to $H$.  A pencil on $X$ is a line $\ell \subset (\mathbb{P}^N)^{\ast}$.  Let $\ell$ be a pencil on $X$ and
  define the variety $$ \widetilde{X}_{\ell} = \left \{ (x,H) \mid x \in \Sigma_H \right \} \subset X \times \ell
  $$ and denote by $\pi: \widetilde{X}_{\ell} \rightarrow \ell \simeq \mathbb{P}^1$ the map given by projection on the $\ell$ factor. Let $$ \begin{array}{ccc} B(\ell) = \bigcap_{H \in \ell} \Sigma_H & ; & S(\ell) = \bigcup_{H \in \ell} Sing(\Sigma_H) \end{array}$$ be the base locus and singular locus of 
  $\ell$ respectively. We have: 
  
  \bigskip
  
  \hspace{-0.6cm} \bf Proposition 5.1: \rm $Sing(\widetilde{X}_{\ell})=(B(\ell) \cap S(\ell)) \times (\ell \cap X^{\ast}) $. 
   
   \bigskip
   
  \hspace{-0.6cm} \bf Proof: \rm Let $(z,t)$ be local
coordinates in a neighborhood $U \times \mathbb{C} \subset X
\times \ell$ around $(x,H_0)$ such that:
$$ \widetilde{X}_{\ell} \cap (U \times \mathbb{C}) = \left \{
(z,t) \mid f_0(z)+tf_1(z)=0 \right
\}$$ where $$ \begin{array}{ccc} \Sigma_{H_0} \cap U = \left \{
f_0 (z)=0 \right \} & , & \Sigma_{H_1} \cap U = \left \{
f_1 (z) =0 \right \} \end{array} $$ with $H_1 \in \ell$
a hyperplane different from $H_0$. The origin is thus a singular point if and only if $f_0(0)=f_1(0)=0$ and 
$df_0(0)=0$. 
$\square$
  
  \bigskip
  
  \hspace{-0.6cm} In particular, since $B(\ell) \subset \Sigma_H$ is an ample divisor, we have: 
  
  \bigskip
  
  \hspace{-0.6cm} \bf Corollary 5.2: \rm $\widetilde{X}_{\ell}$ is singular if $dim(S(\ell)) \geq 1$.
    
  \bigskip
  
    \hspace{-0.6cm} Note that by definition $S(\ell)= \bigcup_{H \in \ell \cap X^{\ast}} Sing(\Sigma_H)$. Recall that a pencil $\ell$ is said to be a Lefschetz 
  pencil on $X$ if it intersects the discriminant variety $X^{\ast}$ transversally. Of course, for such an intersection to exist and be non-void we have to assume that $codim(X^{\ast})=1$. It is well known that if $\ell$ is a Lefschetz pencil 
  the map $\pi: \widetilde{X}_{\ell} \rightarrow \ell$, refered to as a Lefschetz pencil, has only isolated singular points and satisfies the Morse(-Bott) conditions locally around these point. In contrast, when $codim(X^{\ast}) \geq 2 $ a line $ \ell \subset ( \mathbb{P}^N )^{\ast} $ will either be disjoint from $X^{\ast}$, or intersect $X^{\ast}$ non-transversally. In the latter case we obtain $dim(S(\ell)) \geq 1$. 
  
  \hspace{-0.6cm} Thus,  by Corollary 5.2 if $dim(S(\ell)) \geq 1$ the map $\pi : \widetilde{X}_{\ell} \rightarrow \ell$ itself cannot give rise to a Morse-Bott degeneration 
  simply because $\widetilde{X}_{\ell}$ will be singular. One can, however, attempt to resolve the singularities of $\widetilde{X}_{\ell}$.  
    
  \bigskip
  
  \hspace{-0.6cm} \bf Example: \rm Let $z=[z_0:...:z_n]$ be homogenous coordinates on $X=\mathbb{P}^n$ and consider the pencil $\ell$ of quadrics given by $$ \Sigma_{\lambda}=
  \left \{ \lambda_0 F_0(z) + \lambda_1 F_1(z) =0 \right \} \subset \mathbb{P}^n $$ Where $$\begin{array}{ccc} F_0(z)=z^2_0+...+z_n^2 & ; & F_1(z) =-(a_2 z_2^2+a_3 z_3^2 +...+
  a_n z_n^2) \end{array}$$ Where the coefficients $a_i$ are non-zero and $a_i \neq a_j$ for all $2 \leq i,j  \leq n$. In particular,  the quadric $\Sigma_{[0:1]}$ is singular with $$Sing(\Sigma_{[0:1]})=\left \{z_2 = ...=z_n=0 \right \}$$ Moreover, since the rest of the singular 
  fibers are isolated and disjoint of the base lcous we have $$ B(\ell) \cap S(\ell) = B(\ell) \cap Sing(\Sigma_{[0:1]})=\left \{ \begin{array}{c} z_0^2+z_1^2=0 \\ z_2=...=z_n=0 \end{array} \right \}$$ which constitutes of the two points $z_{\pm} = [1,\pm i,0,...,0]$. If one would blow up $\widetilde{X}_{\ell}$ along the two points $\left \{ z_-,z_+ \right \} \times \ell $ one would obtain a degeneration. However, note that the fibres of this degeneration would, at this point, be a quadric blown up at two points and not the original quadrics we began with.

\end{document}